\documentclass[envcountchap, envcountsame]{llncs}
\usepackage{amsmath}
\usepackage{makeidx}  
\usepackage{graphicx}
\usepackage{subfigure}
\usepackage{paralist}
\usepackage{amssymb}
\usepackage {lineno}
\usepackage[utf8]{inputenc}
\usepackage{comment} 
\usepackage{float}

\usepackage{tikz-cd}
\usepackage{tikz}

\title{Resolve the lightcone singularity and extend the $SO(3,1)$-action}

\author{Joachim Wehler}

\institute {LMU Munich (Ludwig-Maximilians-Universit{\"a}t M{\"u}nchen), Germany \\ Department of Mathematics \\ 
Theresienstrasse 39/I, 80333 M{\"u}nchen, Germany\\
\email{Joachim.Wehler@gmx.net} }

\begin{document}
\maketitle
\setcounter{footnote}{0} 

\begin{abstract}
Papadopoulos showed in 2021: The action of the Lorentz group $SO(3,1)$ on the lightcone of Minkowski space lifts to an $SO(3,1)$-action on the blow-up of the lightcone. From this, he draws conclusions regarding further hypothetical vacuum states in quantum field theories.
 
The present note shows that the possbility to lift \mbox{the $SO(3,1)$-action} to the blow-up is a pure \emph{mathematical} result. It does not depend on the context from physics, in particular it does not make use of the Dirac equation. Instead, it is a simple application of the blow-up theory from algebraic geometry. Lifting is possible because the group action is linear and the exceptional divisor of the blow-up is the projectivized tangent cone.
\end{abstract}

\begin{flushleft}
Papadopoulos \cite{Papadopoulos2021} showed that the action of the Lorentz group $SO(3,1)$ on the lightcone $LC$ of Minkowski space lifts to an $SO(3,1)$-action on the \mbox{blow-up $\widetilde{LC}$} of the lightcone. From this, he draws conclusions regarding further hypothetical vacuum states in  quantum field theories.\bigskip

We show that the possbility to lift \mbox{the $SO(3,1)$-action} to the blow-up is a pure \emph{mathematical} result. It does not depend on the context from physics, in particular it does not make use of the Dirac equation. Instead, it is a simple application of the blow-up theory from algebraic geometry. The lifting result relies on the fact that the group action is linear and the exceptional divisor of the blow-up is the projectivized tangent cone.\bigskip 

The light \mbox{cone $LC$}, considered as the boundary of a hollow cone, is \mbox{a $3$-dimensional} real hypersurface of the $4$-dimensional Minkowski \mbox{$\mathbb{R}^{3,1}$}. It is defined as the quadric
$$LC:=\left\{x=(x_0,x_1,x_2,x_3) \in \mathbb{R}^{3,1}:\ x_0^2-x_1^2-x_2^2-x_3^2=0\right\}$$
The lightcone has a unique singularity, namely its vertex, the point 
$$0 \in LC$$ 
Blowing up $LC$ at the singularity constructs a \mbox{non-singular $3$-dimensional} \mbox{variety $\widetilde{LC}$} and a canonical projection
$$\pi: \widetilde{LC} \xrightarrow{} LC$$
with 
$$E:=\pi^{-1}(0) \subset \widetilde{LC}$$
a smooth $2$-dimensional hypersurface, named the \emph{exceptional divisor} of the blow-up. Outside the exceptional divisor the restriction
$$\pi:\widetilde{LC} \setminus E \xrightarrow{} LC \setminus \{0\}$$
is an isomorphism. The points of $E$ correspond bijectively to the \mbox{lines $l \subset LC$} passing through the vertex. Equivalently, the exceptional divisor is the parameter space of all those lines: Each point of the punctured lightcone 
$$x=(x_0,x_1,x_2,x_3) \in LC \setminus \{0\}$$  
satisfies $$x_0\neq 0$$
Hence it determines a unique line $l_x \subset LC$ passing through $x$ and the vertex. The line corresponds to the point of the exceptional divisor with homogeneous coordinates
$$(x_0:x_1:x_2:x_3) \in E \simeq \mathbb{P}(LC),$$
with $\mathbb{P}(LC)$ the projectivization of the cone $LC$. The last isomorphy is a special case of the isomorphy  
\begin{itemize}
\item between the exceptional divisor of the blow-up of a variety at a given point \bigskip

\item and the projectivized tangent cone of the variety at this point;
\end{itemize}
e.g., see \cite[Exercise IV-24]{EH2000}. The projective variety $\mathbb{P}(LC)$ is the quadric hypersurface in $\mathbb{P}^3$
$$\mathbb{P}(LC) = \{(x_0:x_1:x_2:x_3) \in \mathbb{P}^3:\ x_0^2-x_1^2-x_2^2-x_3^2=0\}$$ 
The Lorentz group $SO(3,1)$ acts on the lightcone $LC$ with fixed point the vertex. The action is linear: If a Lorentz \mbox{transformation $S\in SO(3,1)$} maps the point $x \in LC \setminus \{0\}$ to the \mbox{point $x^\prime \in LC$}, then also $x^\prime \neq 0$, and  $S$ maps the whole line $l_x$ isomorphically onto the line $l_{x^\prime}$. Therefore $SO(3,1)$ also acts on the set of all lines of $LC$ passing through the vertex. They make up the tangent cone \mbox{of $LC$} at the vertex. Each line $l_x$ is determined by the homogeneous \mbox{coordinates $(x_0:x_1:x_2:x_3)$} of $x$. Hence \mbox{the $SO(3,1)$-action} \mbox{on $LC$} extends to a smooth action on $\widetilde{LC}$ 
$$SO(3,1) \times{} \widetilde{LC} \xrightarrow{} \widetilde{LC}$$
when defining the action on the exceptional divisor as 
$$SO(3,1) \times{} \mathbb{P}(LC) \xrightarrow{} \mathbb{P}(LC),\ (S,\ (x_0:x_1:x_2:x_3)) \mapsto (x_0^\prime:x^\prime_1:x^\prime_2:x^\prime_3)$$
The following diagram commutes 
$$\begin{tikzpicture}
  \matrix (m) [matrix of math nodes,row sep=4em,column sep=4em,minimum width=3em] {
    	SO(3,1) \times{} \widetilde{LC} 	&  \widetilde{LC}	 				\\
	SO(3,1) \times{} LC				&			LC						\\	} ;
   \path[-stealth]
    (m-1-1) edge 			node 		[above] 		{} 						(m-1-2)
    (m-1-1) edge 			node 		[left]     		{$id \times \pi$} 		(m-2-1) 
    (m-1-2) edge 			node 		[right]     		{$\pi$} 				(m-2-2) 
    (m-2-1) edge		 	node 		[above]   		{} 						(m-2-2) ;  
\end{tikzpicture}$$
and the exceptional divisor 
$$\pi^{-1}(0)=E \simeq \mathbb{P}(LC)$$ 
is an additional orbit of the $SO(3,1)$-action.
\end{flushleft}

\bibliographystyle{splncs03} 
{}


\end{document}